\documentclass{article}

\usepackage{amsmath,amsfonts,amsthm,amssymb,graphicx}
\usepackage[all,2cell,ps]{xy}

\theoremstyle{plain}
\newtheorem{thm}{Theorem}[section]
\newtheorem{lem}[thm]{Lemma}
\newtheorem{prop}[thm]{Proposition}
\newtheorem{cor}[thm]{Corollary}

\newtheorem{conjec}[thm]{Conjecture}

\theoremstyle{definition}

\newtheorem*{rem}{Remark}
\newtheorem*{rems}{Remarks}

\theoremstyle{remark}

\newtheorem{clm}{Claim}

\newcommand{\al}{\alpha}
\newcommand{\F}{\mathbb F}
\newcommand{\Z}{\mathbb Z}

\newcommand{\R}{\mathbb R}

\newcommand{\Q}{\mathbb Q}
\newcommand{\B}{\mathbb B}

\newcommand{\de}{\delta}

\newcommand{\Sig}{\Sigma}

\newcommand{\Lam}{\Lambda}
\newcommand{\A}{\mathbb{A}}
\newcommand{\C}{\mathbb C}

\newcommand{\Gam}{\Gamma}

\newcommand{\conj}{\overline}
\newcommand{\om}{\omega}

\newcommand{\Del}{\Delta}
\newcommand{\mc}{\mathcal}
\newcommand{\del}{\delta}

\DeclareMathOperator{\PSL}{PSL}

\DeclareMathOperator{\PGL}{PGL}

\DeclareMathOperator{\PU}{PU}
\DeclareMathOperator{\SU}{SU}
\DeclareMathOperator{\PSU}{PSU}

\DeclareMathOperator{\para}{P}

\newcommand{\mf}{\mathfrak}

\DeclareMathOperator{\Aut}{Aut}

\newcommand{\ssm}{\smallsetminus}

\DeclareMathOperator{\Disc}{D}

\newcommand{\wt}{\widetilde}

\newenvironment{pf}{\begin{proof}}{\end{proof}}

\newenvironment{enum}{\begin{enumerate}}{\end{enumerate}}

\title{Hurwitz ball quotients}
\author{Matthew Stover\footnote{This material is based upon work supported by the National Science Foundation 
under Grant Number NSF 0943832.} \\ \small{Temple University}\\ \small{\textsf{mstover@temple.edu}}}
\date{\today}

\begin{document}

\maketitle

\begin{abstract}
We consider the analogue of Hurwitz curves, smooth projective curves $C$ of genus $g \ge 2$ that realize equality in the Hurwitz bound $|\Aut(C)| \le 84 (g - 1)$, to smooth compact quotients $S$ of the unit ball in $\C^2$. When $S$ is arithmetic, we show that $|\Aut(S)| \le 288 e(S)$, where $e(S)$ is the (topological) Euler characteristic, and in the case of equality show that $S$ is a regular cover of a particular Deligne--Mostow orbifold. We conjecture that this inequality holds independent of arithmeticity, and note that work of Xiao makes progress on this conjecture and implies the best-known lower bound for the volume of a complex hyperbolic $2$-orbifold.
\end{abstract}


\section{Introduction}\label{sec:intro}


The classical Hurwitz bound for the order of the automorphism group of a smooth projective curve $C$ of genus $g \geq 2$ is
\[
|\Aut(C)| \leq 84(g - 1) = - 42 e(C),
\]
where $e$ denotes the (topological) Euler characteristic. Curves realizing this bound are called \emph{Hurwitz curves}, and their fundamental groups are the normal torsion-free finite index subgroups of the $(2, 3, 7)$-triangle group, $\Del_{2,3,7}$. The most famous Hurwitz curve is the Klein quartic, which has automorphism group $\PSL_2(\F_7)$ and fundamental group the congruence subgroup of level $7$ in the arithmetic Fuchsian group $\Del_{2,3,7}$. One also has that \emph{Hurwitz groups}, the automorphism groups of Hurwitz curves, are in one-to-one correspondence with finite groups generated by two elements $x$ and $y$ such that $x$ has order $2$, $y$ has order $3$, and $x y$ has order $7$.

The purpose of this paper is to study the analogous problem for smooth projective surfaces $S$ uniformized by the complex $2$-ball $\B^2$ with its Bergman metric, that is, complex hyperbolic $2$-space. In other words, we study how large the automorphism group of a ball quotient can be relative to its Euler characteristic. When $\pi_1(S) \subset \PU(2, 1)$ is arithmetic, we prove the following.


\begin{thm}\label{thm:ArithmeticAut}
Let $S$ be a smooth projective surface uniformized by $\B^2$, and suppose that $\pi_1(S) \subset \PU(2, 1)$ is an arithmetic lattice. Then $|\Aut(S)| \le 288 e(S)$ with equality if and only if $S / \Aut(S)$ is the Deligne--Mostow orbifold associated with the ball tuple $(\frac{2}{12}, \frac{2}{12}, \frac{2}{12}, \frac{7}{12}, \frac{11}{12})$.
\end{thm}


See \S\ref{sec:DM} for more on this Deligne--Mostow orbifold, which is entry $63$ in the appendix to \cite{Mostow} and entry $61$ in \cite{Thurston}. We note that this orbifold also appears in Mostow's earlier work on lattices generated by complex reflections  \cite{MostowReflect}, and is associated with the complex reflection group he calls $\Gam(3, \frac{1}{3})$ (see \cite[\S 7]{Sauter}). Both interpretations make it a natural generalization of the $(2, 3, 7)$-triangle orbifold. As we will see, Theorem \ref{thm:ArithmeticAut} is equivalent to the following.


\begin{thm}\label{thm:DMVolume}
Let $\mc{O} = \B^2 / \Gamma$ be a finite volume quotient of $\B^2$ by an arithmetic lattice $\Gamma \subset \PU(2, 1)$. Then
\[
e(\mc{O}) \ge \frac{1}{288}
\]
with equality if and only if $\mc{O}$ is the Deligne--Mostow orbifold associated with the ball tuple $(\frac{2}{12}, \frac{2}{12}, \frac{2}{12}, \frac{7}{12}, \frac{11}{12})$.
\end{thm}


Here $e(\mc{O})$ is the orbifold Euler characteristic. We conjecture that Theorems \ref{thm:ArithmeticAut} and \ref{thm:DMVolume} hold without the arithmetic assumption. We find this natural given that the minimum volume hyperbolic $2$- and $3$-orbifolds are both arithmetic \cite{Gehring--MartinMRL, Gehring--Martin, Marshall--Martin}. In particular, we conjecture the following.


\begin{conjec}\label{MinVolConjecture}
The Deligne--Mostow orbifold associated with the ball tuple $(\frac{2}{12}, \frac{2}{12}, \frac{2}{12}, \frac{7}{12}, \frac{11}{12})$ is the unique minimum volume ball quotient.
\end{conjec}


See \S \ref{sec:AutGps} for more on this conjecture and an equivalent conjecture analogous to Theorem \ref{thm:ArithmeticAut} stated in terms of automorphism groups. We also note that work of Xiao implies the best-known unconditional lower bound for the volume of a closed complex hyperbolic $2$-orbifold, namely $\pi^2 / 1944 \approx 0.005077$; see Proposition \ref{prop:VolLowerBound}.

We prove Theorems \ref{thm:ArithmeticAut} and \ref{thm:DMVolume} in three steps.
\begin{enum}

\item In \S \ref{sec:AutGps}, we explain the relationship between automorphism groups of manifolds and topological invariants, in particular the equivalence of Theorem \ref{thm:ArithmeticAut} and Theorem \ref{thm:DMVolume}. We also discuss Conjecture \ref{MinVolConjecture} further and give an equivalent restatement in terms of automorphism groups.

\item In \S \ref{sec:Arithmetic}, we prove that there is a unique lattice $\Gam$ of minimal covolume amongst all arithmetic lattices in $\PU(2, 1)$. The proof, which relies on Prasad's volume formula \cite{Prasad}, can be pieced together from results contained in Prasad and Yeung's paper on fake projective planes \cite{Prasad--Yeung}. We give a complete streamlined proof that avoids some technical points in their paper.

\item Finally, in \S \ref{sec:DM}, we describe Deligne and Mostow's orbifolds and prove that the lattice $\Gam$ from Step $2$ is the orbifold fundamental group of the Deligne--Mostow orbifold with ball tuple $(\frac{2}{12}, \frac{2}{12}, \frac{2}{12}, \frac{7}{12}, \frac{11}{12})$. This completes the proofs of Theorems \ref{thm:ArithmeticAut} and \ref{thm:DMVolume}.

\end{enum}

In \S \ref{sec:Gps}, we study the automorphism groups that are extremal for Theorem \ref{thm:ArithmeticAut}. Recall that Hurwitz groups, the automorphism groups of curves of genus $g \ge 2$ with exactly $84(g - 1)$ automorphisms, are in one-to-one correspondence with groups generated by two elements $x, y$ such that $x$ has order $2$, $y$ has order $3$, and $x y$ has order $7$. The situation for arithmetic quotients of $\B^2$ is similar, with the extremal automorphism groups related to groups generated by the exceptional group of order $288$ in the Shephard--Todd classification of finite complex reflection groups and an element of order $3$. A precise statement is complicated, so we refer the interested reader to \S \ref{sec:Gps} for details.

We also explore the smallest surfaces that are extremal for Theorem \ref{thm:ArithmeticAut}. For example, we show the following.


\begin{thm}\label{thm:Smallest}
The smooth arithmetic ball quotient $S$ with $|\Aut(S)| = 288 e(S)$ and minimal Euler characteristic amongst all such examples is a quotient of the ball by a principal congruence lattice. It has Euler characteristic $63$ and automorphism group $\PSU(3, \F_3) \times (\Z / 3 \Z)$.
\end{thm}


See \S \ref{sec:Gps} for more about this surface, e.g., its numerical invariants. This surface is a quotient of the ball by a principal congruence subgroup in the commensurability class of arithmetic subgroups of $\PU(2, 1)$ defined by a hermitian form on $\Q(\zeta)^3$, where $\zeta$ is a primitive $12^{th}$ root of unity (i.e., the commensurability class of the Deligne--Mostow orbifold in Theorem \ref{thm:ArithmeticAut}). One can also realize it as a connected component of a PEL Shimura variety. This gives a number of analogies to the Klein quartic. We also give a description of the second smallest such surface, which has automorphism group $\PSU(3, \F_3) \times \mathrm{A}_4$, where $\mathrm{A}_4$ is the alternating group on four letters, and is a $4$-fold cover of the surface in Theorem \ref{thm:Smallest}.


\subsubsection*{Acknowledgments} I am indebted to Domingo Toledo for discussions about visualizing the orbifold structure on a Deligne--Mostow quotient. Any insight into the geometry of these spaces not in the standard literature should be considered his, not mine. I also want to thank the referee for suggestions that undoubtedly improved this paper.


\section{Automorphism groups and topological invariants}\label{sec:AutGps}


Let $M$ be a closed Riemannian manifold of negative sectional curvature and $\wt{M}$ be its universal cover. Then $\Aut(M)$ is well-known to be finite and $\mc{O} = M / \Aut(M)$ is a Riemannian orbifold with universal cover $\wt{M}$, that is, there is a discrete cocompact group $\Gam$ of isometries of $\wt{M}$ such that $\mc{O} = \wt{M} / \Gam$. The following simple lemma is key to our approach to bounding orders of automorphism groups.


\begin{lem}\label{lem:VolAutBound}
Let $\wt{M}$ be a complete simply connected Riemannian manifold with negative sectional curvature, and suppose that there is a constant $c > 0$ such that $\mathrm{vol}(\wt{M} / \Gam) \ge c$ for every discrete cocompact group of isometries $\Gam$ of $\wt{M}$. Then $|\Aut(M)| \le \mathrm{vol}(M) / c$ for every compact Riemannian manifold $M$ with universal cover $\wt{M}$.
\end{lem}

\begin{pf}
The lemma follows immediately from the following calculation:
\[
c \le \mathrm{vol}(\mc{O}) = \frac{1}{|\Aut(M)|} \mathrm{vol}(M),
\]
where $\mc{O} = M / \Aut(M)$.
\end{pf}


By the Margulis lemma, Lemma \ref{lem:VolAutBound} applies to every symmetric space of noncompact type and, more generally, complete manifolds of bounded negative curvature (see \cite{Kazhdan--Margulis} and \cite[\S 10]{BGS}). Clearly an analogous statement holds for any topological invariant proportional to volume. This gives the following consequence of Lemma \ref{lem:VolAutBound}.


\begin{cor}\label{cor:BallInvariantAut}
Let $\Gam < \PU(2, 1)$ be a torsion-free cocompact lattice and $S = \B^2 / \Gam$ be the associated smooth projective surface. Then there is a constant $b > 0$, independent of $S$, such that
\begin{equation}\label{eq:BallInvariantAut}
|\Aut(S)| \le b e(S) = b c_2(S).
\end{equation}
\end{cor}

\begin{pf}
Consider the metric on $S$ of constant holomorphic curvature $-1$. Chern--Gauss--Bonnet implies that
\[
\mathrm{vol}(S) = \frac{8 \pi^2}{3} e(S),
\]
so we take $b = 8 \pi^2 / 3 c$, where $c$ is the constant in Lemma \ref{lem:VolAutBound}.
\end{pf}


We call the extremal surfaces for this inequality \emph{Hurwitz ball quotients}. Theorem \ref{thm:ArithmeticAut} says that $b = 288$ suffices when $S$ is arithmetic, hence we can take $c = \pi^2 / 108 \approx 0.0914$ in Lemma \ref{lem:VolAutBound}. Conjecture \ref{MinVolConjecture} is equivalent to the conjecture that $b = 288$ is optimal without the arithmetic assumption. More specifically, we conjecture the following, which is equivalent to Conjecture \ref{MinVolConjecture}.


\begin{conjec}\label{AutConjecture}
Let $S$ be a smooth projective surface uniformized by $\B^2$. Then
\begin{equation}\label{eq:AutConjecture}
|\Aut(S)| \le 288 e(S).\tag{$\diamondsuit$}
\end{equation}
Furthermore, the following statements are equivalent:
\begin{enum}

\item The surface $S$ is extremal for \eqref{eq:AutConjecture}, i.e., $|\Aut(S)| = 288 e(S)$.

\item The complex hyperbolic $2$-orbifold $S / \Aut(S)$ is the Deligne--Mostow orbifold associated with the ball tuple $(\frac{2}{12}, \frac{2}{12}, \frac{2}{12}, \frac{7}{12}, \frac{11}{12})$.

\end{enum}
Equivalently, the minimal volume of a compact complex hyperbolic $2$-manifold is $\pi^2 / 108$, and this volume is uniquely realized by the above orbifold.
\end{conjec}


\begin{rem}
Using standard numerical arguments in algebraic geometry, we can rephrase all the above results in terms of either the holomorphic Euler characteristic $\chi(\mc{O}_S)$ or the self-intersection of the canonical bundle $K_S^2 = c_1^2(S)$. For ball quotients, we have the following relationships between these invariants and $e(S)$:
\begin{eqnarray}
\chi(\mc{O}_S) &=& \frac{1}{3} e(S) \label{eq:HolECProp} \\
K_S^2 &=& 3 e(S) \label{eq:K2ECProp}
\end{eqnarray}
To prove these, Hirzebruch proportionality for $\B^2$ gives $c_1^2(S) = 3 c_2(S)$, which is exactly \eqref{eq:K2ECProp}, and \eqref{eq:HolECProp} follows from \eqref{eq:K2ECProp} and Noether's formula \cite[p.~432]{Griffiths--Harris}.
\end{rem}


For any smooth surface of general type, Xiao \cite{Xiao, Xiao2} showed that
\[
|\Aut(S)| \leq \begin{cases} 42^2 c_1^2(S) & \widehat{S / \Aut(S)}\ \textrm{rational} \\ \\ 288 c_1^2(S) & \textrm{otherwise}, \end{cases}
\]
where $\widehat{\ }$ denotes the resolution of singularities of the possibly singular variety $S / \Aut(S)$. This is a natural generalization of the Hurwitz bound, since $c_1(C) = \chi(\mc{O}_C)$ for curves, so the Hurwitz bound is $|\Aut(C)| \le 42 |c_1(C)|$. Further, the assumption that $S$ has general type corresponds with the assumption $g \geq 2$ in dimension one. Xiao also proves that $|\Aut(S)| = 42^2 c_1^2(S)$ if and only if $S$ is a quotient of $C \times C$ with $C$ a Hurwitz curve.

Since a product of curves is not uniformized by $\B^2$, arguments in \cite{Xiao2} (see \S 10) show that we can actually replace $42^2$ with $1728$. Applying Hirzebruch proportionality and Chern--Gauss--Bonnet, Xiao's results have the following consequence.


\begin{prop}\label{prop:VolLowerBound}
Let $\mc{O}$ be a closed complex hyperbolic $2$-orbifold. Then
\[
\mathrm{vol}(\mc{O}) \ge \frac{\pi^2}{1944} \approx 0.005077
\]
in the metric of constant holomorphic curvature $-1$. Equivalently, we can take $b = 5184$ in Corollary \ref{cor:BallInvariantAut}.
\end{prop}


The best previous unconditional lower bound for the volume of a complex hyperbolic $2$-orbifold is $2.918 \times 10^{-9}$ from Adeboye--Wei \cite{Adeboye--Wei}. For the noncompact finite volume case, see \cite{Parker, StoverVol}. Recall that Conjecture \ref{MinVolConjecture} implies that the minimal volume of a complex hyperbolic $2$-orbifold is exactly $\pi^2 / 108 \approx 0.0913853$, and that this volume is uniquely realized by the Deligne--Mostow orbifold in the conjecture.


\section{Arithmetic lattices of minimal covolume}\label{sec:Arithmetic}


In this section we determine the cocompact arithmetic lattice $\Gam < \PU(2, 1)$ of minimal covolume. Since related results appear in many places in the literature (e.g., \cite{Prasad--Yeung, Emery--Stover}), we cut straight to the chase.


\begin{thm}\label{thm:MinVol}
There is a unique arithmetic lattice $\Gam < \PU(2, 1)$ such that
\[
e(\B^2 / \Lam) \ge e(\B^2 / \Gam) = \frac{1}{288}
\]
for every arithmetic lattice $\Lam < \PU(2, 1)$ with equality if and only if $\Lam \cong \Gam$, where $e$ denotes the orbifold Euler characteristic. The commensurability class of $\Gam$ (that is, the associated $\Q$-algebraic group) is uniquely determined by the following data:
\begin{enum}

\item $k = \Q(\al)$, where $\al^2 = 3$

\item $\ell = k(\beta)$, where $\beta^2 = -1$

\item $\tau$ the nontrivial Galois involution of the quadratic extension $\ell / k$

\item $h$ the $\tau$-hermitian form on $\ell^3$ with matrix
\[
\begin{pmatrix}
1 & 0 & 0 \\ 0 & 1 & 0 \\ 0 & 0 & 1 - \al
\end{pmatrix}
\]

\end{enum}
\end{thm}

\begin{pf}
Let $\Gam \subset \PU(2, 1)$ be an arithmetic lattice. Associated with $\Gam$ there is a totally real field $k$ of degree $n$ over $\Q$ and a simple simply connected $k$-algebraic group $G$ such that $\mathrm{Res}_{k / \Q}(G)(\R) \cong \SU(2, 1) \times \SU(3)^{n - 1}$, where $\mathrm{Res}_{k / \Q}$ denotes Weil restriction of scalars. Moreover, there is a totally imaginary quadratic extension $\ell$ of $k$, a central simple $\ell$-division algebra $D$ with involution $\tau$ of second kind (i.e., such that $\tau$ restricts to the Galois involution of $\ell / k$) and degree $d \in \{1, 3\}$, and a nondegenerate $\tau$-hermitian form $h$ on $D^r$, where $d r = 3$, such that $G$ is the special unitary group of $h$. For example, when $D$ has degree $1$, $G$ is the special unitary group of a hermitian form $h$ on $\ell^3$ such that $h$ is indefinite at exactly one place of $\ell$ (i.e., complex conjugate pair of embeddings $\ell \to \C$). It is not hard to show that $G$ is uniquely determined up to $k$-isomorphism by $D$, that is, that the choice of involution and hermitian form are irrelevant; see \cite[\S 1.2]{Prasad--Yeung}.

Denote by $V = V_\infty \cup V_f$ (resp.~$W = W_\infty \cup W_f$) the places of $k$ (resp.~$\ell$), decomposed into its subsets of infinite and finite places. For $v \in V$ (resp.~$w \in W$), $k_v$ (resp.~$\ell_w$) denotes the completion of $k$ (resp.~$\ell$) at $v$ (resp.~$w$). For $v \in V_f$, let $q_v$ be the order of the residue field of $k_v$. The absolute value of the discriminant of a field $F$ is denoted $\Disc_F$, the Dedekind $\zeta$-function of a field $F$ is $\zeta_F(s)$, the $L$-function of an extension $E / F$ is $L_{E / F}(s)$, $\A_F$ is the adele ring of a field $F$, and $\A_{F, f}$ the subring of finite adeles.

Let
\[
\para = \prod_{v \in V_f} \para_v \subset G(\A_{k, f})
\]
be a coherent family of parahoric subgroups and $\Gam_{\para} = \para \cap G(k)$ be the associated lattice. See \cite{Prasad} for the precise definition of a coherent family of parahoric subgroups. The projection of $\Gam_{\para}$ to $\PU(2, 1)$ is a lattice commensurable with our original lattice $\Gam$. Following \cite[\S 2.4]{Prasad--Yeung}, we now give the formula for the orbifold Euler characteristic $e(\B^2 / \Gam_{\para})$ via Prasad's volume formula \cite{Prasad}.

As in \cite[\S 2.2]{Prasad--Yeung}, let $\mc{T}$ be the set of places $v \in V_f$ such that (1) $v$ is unramified in $\ell$ and (2) $\para_v$ is not a hyperspecial parahoric subgroup of $G(k_v)$. Define constants $e(\para_v)$ and $e^\prime(\para_v)$ as in \cite{Prasad--Yeung}. We do not need these constants in what follows, only the facts that $e(\para_v), e^\prime(\para_v) \ge 1$ and equal $1$ when $\para_v$ is maximal hyperspecial, so we leave it to the interested reader to compute them. Applying Prasad's volume formula \cite[\S 3.7]{Prasad}, we see that
\begin{equation}\label{eq:PrincipalVol}
e(\B^2 / \Gam_{\para}) = \frac{9 \Disc_\ell^{5/2} \zeta_k(2) L_{\ell / k}(3)}{(16 \pi^5)^n \Disc_k} \prod_{v \in \mc{T}} e^\prime(\para_v),
\end{equation}
where, as always, $e(X)$ is the orbifold Euler characteristic of $X$. See the discussion in \cite[\S 1.3]{Prasad--Yeung}, and note that we are considering the quotient of $\B^2$ not $\SU(2, 1)$ as there.

It is known that the lift $\wt{\Gam}$ of the maximal arithmetic lattice $\Gam$ to $\SU(2, 1)$ is the normalizer in $\SU(2, 1)$ of $\Gam_{\para}$ for some coherent family of parahoric subgroups $\para$. Let $\mc{T}_0 \subset \mc{T}$ be the subset of places $v$ where $D_v = D \otimes_k k_v$ is still a division algebra. This equals the set of places where $G(k_v)$ is anisotropic, and such a place $v$ necessarily splits in $\ell$. By \cite[\S 2.3]{Prasad--Yeung}, we have
\begin{equation}\label{eq:IndexBound}
[\wt{\Gam} : \Gam_{\para}] \le 3^{1 + \#\mc{T}_0} h_{\ell, 3} \prod_{v \in \mc{T} \ssm \mc{T}_0} \#\Xi_{\Theta_v},
\end{equation}
where $h_{\ell, 3}$ is the order of the $3$-primary part of the class group of $\ell$, and $\#\Xi_{\Theta_v} = 1$ unless $v$ splits in $\ell$ and $\para_v$ is an Iwahori subgroup of $G(k_v)$, in which case it equals $3$ (see \cite[\S 2.2]{Prasad--Yeung}). We then define:
\[
e^{\prime \prime}(\para_v) = \begin{cases} e^\prime(\para_v) / 3 & v \in \mc{T}_0 \\ \\ e^\prime(\para_v) / \# \Xi_{\Theta_v} & v \in \mc{T} \ssm \mc{T}_0 \end{cases}
\]

Also, by \cite[\S 3.10]{Prasad}, $e^{\prime \prime}(\para_v) \ge 1$ for all $v$. Therefore, for any $v \notin \mc{T}_0$ we could have chosen $\para_v$ hyperspecial and the resulting lattice would be commensurable $\Gam_{\para}$ and have covolume less than or equal to that of $\Gam_{\para}$. Consequently, we assume from here forward that
\begin{enum}

\item $\mc{T} = \mc{T}_0$ is the set of places $v$, necessarily split in $\ell$, where $G(k_v)$ is anisotropic,

\item $\para_v$ is maximal for all $v$ that ramify in $\ell$, and

\item $\para_v$ is maximal and hyperspecial for all remaining $v \in V_f$.

\end{enum}
Let $\wt{\Gam} \subset \SU(2, 1)$ be the maximal lattice that normalizes $\Gam_{\para}$. Combining \eqref{eq:PrincipalVol} and \eqref{eq:IndexBound} implies that
\begin{equation}\label{eq:VolBound}
e(\B^2 / \wt{\Gam}) \ge \frac{3 \Disc_\ell^{5/2} \zeta_k(2) L_{\ell / k}(3)}{(16 \pi^5)^n \Disc_k h_{\ell, 3}} \prod_{v \in \mc{T}} e^{\prime \prime}(\para_v).
\end{equation}

Now, suppose that $\Lam \subset \PU(2, 1)$ is an arithmetic lattice such that $e(\B^2 / \Lam) \le 1 / 288$. We want to show that there is a unique such $\Lam$, and that it comes from the construction in the statement of the theorem. Suppose that $\Gam_{\para}$ is commensurable with the lift of $\Lam$ to $\SU(2, 1)$ and that $\para$ satisfies assumptions (1)-(3) above. If $\wt{\Gam}$ is the normalizer of $\Gam_{\para}$ in $\SU(2, 1)$, then $\wt{\Gam}$ has minimal covolume amongst all lattices in its commensurability class by our assumptions on $\para$. Therefore, we also have $e(\B^2 / \wt{\Gam}) \le 1 / 288$. Applying \eqref{eq:VolBound}, we have
\begin{equation}\label{eq:288VolBound}
\frac{\Disc_\ell^{5/2} \zeta_k(2) L_{\ell / k}(3)}{(16 \pi^5)^n \Disc_k h_{\ell, 3}} \prod_{v \in \mc{T}} e^{\prime \prime}(\para_v) \le \frac{1}{864}.
\end{equation}
Since $e^{\prime \prime}(\para_v) \ge 1$, \eqref{eq:288VolBound} implies that
\begin{equation}\label{eq:FieldBound}
\frac{\Disc_\ell^{5/2} \zeta_k(2) L_{\ell / k}(3)}{(16 \pi^5)^n \Disc_k h_{\ell, 3}} \le \frac{1}{864}.
\end{equation}

It is a well-known consequence of the Brauer--Siegel Theorem that
\begin{equation}\label{eq:BrauerSiegel}
h_{\ell, 3} \le h_\ell \le \frac{w_\ell}{R_\ell} s (s - 1) \Gam(s)^n \left( \frac{\Disc_\ell}{(2 \pi)^{2 n}} \right)^{s / 2} \zeta_\ell(s)
\end{equation}
for every $s > 1$, where $h_\ell$ is the class number of $\ell$, $R_\ell$ is the regulator, $w_\ell$ is the order of the group of roots of unity in $\ell$, and $\Gam(s)$ is the usual gamma function. Applying Slavutskii's bound $R_\ell \ge 0.00136 e^{0.57 n} w_\ell$ \cite{Slavutskii}, we have
\begin{equation}\label{eq:BrauerSiegel2}
h_{\ell, 3} \le h_\ell \le \frac{s (s - 1) \Gam(s)^n}{0.00136 e^{0.57 n}} \left( \frac{\Disc_\ell}{(2 \pi)^{2 n}} \right)^{s / 2} \zeta_\ell(s).
\end{equation}
Factoring all this into \eqref{eq:FieldBound}, we get
\begin{equation}\label{eq:NearBigBound}
\frac{(0.00136)e^{0.57 n} 2^{(s - 4)n} \pi^{(s - 5) n} \Disc_\ell^{(5 - s) / 2} \zeta_k(2) L_{\ell / k}(3)}{s (s - 1) \Gam(s)^n \zeta_\ell(s) \Disc_k} \le \frac{1}{864}
\end{equation}
for all $s > 1$.

Further, $\zeta_k(2) L_{\ell / k}(3) > \zeta(2 n)^{1 / 2}$ by \cite[Cor.~2.8]{Prasad--Yeung} and $\zeta_\ell(s) \le \zeta(s)^n$. Applying these facts and rearranging as in \cite[\S 2.8]{Prasad--Yeung}, we get
\begin{equation}\label{eq:BigBound}
\Disc_\ell^{1 / 2 n} \le 2 \left( \Gam(1 + \de) \zeta(1 + \de) \pi^{4 - \de} e^{-0.57} \right)^{\frac{1}{3 - \del}} \left( \frac{\de (\de + 1)}{(1.17504) \zeta(2 n)^{1 / 2}} \right)^{\frac{1}{(3 - \de) n}}
\end{equation}
for any $0 < \de \le 2$. Note that the right-hand side of \eqref{eq:BigBound} also bounds $\Disc_k^{1/n}$, since $\Disc_\ell^2 / \Disc_k$ is always an integer.

\begin{clm}
If \eqref{eq:NearBigBound} holds, then $n = [k : \Q] \le 2$. That is, $k$ is either $\Q$ or real quadratic.
\end{clm}

\begin{pf}
In \cite[\S 7]{Prasad--Yeung}, Prasad and Yeung show that if $[k : \Q] \ge 6$, then no $k$-algebraic group as above can produce a space with orbifold Euler characteristic less than or equal to $3$. In particular, it cannot give one with orbifold Euler characteristic $1 / 288$. Therefore, we focus on ruling out the cases $n = 3, 4, 5$ using \eqref{eq:BigBound}. The motivated reader can also mimic the argument given below to rule out $n \ge 6$.

For $1 \le n \le 5$, Table \ref{tb:DelBound} gives a particularly nice choice of $\de$ for \eqref{eq:BigBound} and the associated bound. We include $n = 1, 2$, since we need those bounds later.
\begin{table}[h]
\begin{center}
\begin{tabular}{|c|c|c|}
\hline
$n$ & $\de$ & $\Disc_\ell^{1 / 2 n} \le$ \\
\hline
1 & $0.00145$ & $6.64809$ \\
2 & $0.395731$ & $9.96044$ \\
3 & $0.523748$ & $10.404$ \\
4 & $0.589587$ & $10.523$ \\
5 & $0.629827$ & $10.5646$ \\
\hline
\end{tabular}
\end{center}
\caption{Choices of $\de$ and the associated upper bound}\label{tb:DelBound}
\end{table}
Recall that the Hilbert class field of $\ell$ is an abelian Galois extension $h(\ell)$ such that $h(\ell) / \ell$ is totally unramified. It is well-known that $[h(\ell) : \ell] = h_\ell$, the class number, and $\Disc_{h(\ell)} = \Disc_\ell^{h_\ell}$. Thus we can also use Table \ref{tb:DelBound} to bound $h_\ell$, and thus $h_{\ell, 3}$. Note that $h(\ell)$ is totally complex. Using the tables of upper bounds for $\Disc_F^{1/[F : \Q]}$ given by Diaz y Diaz \cite{DiazyDiaz}, which are derived from the well-known Odlyzko bound \cite{Odlyzko}, we get the bounds found in Table \ref{tb:ClassBound}.
\begin{table}[h]
\begin{center}
\begin{tabular}{|c|c|}
\hline
$n$ & $h_{\ell, 3} \le$ \\
\hline
1 & $3$ \\
2 & $3$ \\
3 & $3$ \\
4 & $1$ \\
5 & $1$ \\
\hline
\end{tabular}
\end{center}
\caption{Class number upper bounds}\label{tb:ClassBound}
\end{table}

Now, for each $n = 3, 4, 5$ we revisit \eqref{eq:FieldBound} with our improved bound on $h_{\ell, 3}$, and obtain:
\begin{equation}\label{eq:FieldBound2}
\frac{\Disc_\ell^{5/2} \zeta_k(2) L_{\ell / k}(3)}{\Disc_k} \le \begin{cases} \frac{(16 \pi^5)^3}{288} & n = 3 \\ & \\ \frac{(16 \pi^5)^n}{864} & n = 4, 5 \end{cases}
\end{equation}
Since $\zeta_k(2) L_{\ell / k}(3) \ge \zeta(2 n)^{1/2}$ and $\Disc_\ell^{1/2} / \Disc_k \ge 1$, we have:
\begin{equation}\label{eq:FieldBound3}
\Disc_\ell^{1/2 n} \le \begin{cases} \frac{(16 \pi^5)^{1/4}}{\zeta(6)^{1/24} 288^{1/12}} \approx 5.214 & n = 3 \\ & \\ \frac{(16 \pi^5)^{1/4}}{\zeta(8)^{1/32} 864^{1/16}} \approx 5.481 & n = 4 \\ & \\  \frac{(16 \pi^5)^{1/4}}{\zeta(10)^{1/40} 864^{1/20}} \approx 5.965 & n = 5 \end{cases}
\end{equation}
The tables in \cite{DiazyDiaz} then rule out $n = 4, 5$, since $\Disc_\ell^{1/2n}$ is bounded below by $5.659$ and $6.600$, respectively. Since $\Disc_k^{1/n} \le \Disc_\ell^{1/2 n}$, we see that if $n = 3$, then $k$ is a totally real cubic field of absolute discriminant at most $141$ and $\ell$ an imaginary quadratic extension of absolute discriminant at most $20102$.

There are exactly two such real cubic fields, $k_1$ with discriminant $49$ and $k_2$ with discriminant $81$, and twelve possible totally complex sextic fields. See \cite{Tables}. All these fields have class number one. Also, $\ell$ must be a quadratic extension of $k$, which leaves us with only
\[
(\Disc_k, \Disc_\ell) \in \left\{ (49, 16807), (81, 19683) \right\}.
\]
We then consider \eqref{eq:FieldBound} and see that
\[
\frac{\Disc_\ell^{5/2} \zeta_k(2) L_{\ell / k}(3)}{(16 \pi^5)^3 \Disc_k} \ge \frac{\Disc_\ell^{5/2} \zeta(6)^{1/2}}{(16 \pi^5)^3 \Disc_k} \approx \begin{cases} 0.00642 & \Disc_k = 49 \\ & \\ 0.00577 & \Disc_k = 81 \end{cases}
\]
is a lower bound for the minimal orbifold Euler characteristic of an arithmetic lattice coming from a $k$-algebraic group as above. These lower bounds are both greater than $1/864$, so neither case can lead to an orbifold of Euler characteristic $1/288$. This rules out $n = 3$, and completes the proof of the claim.
\end{pf}

\begin{clm}
We cannot have $k = \Q$.
\end{clm}

\begin{pf}
In the case $n = 1$, we have $h_{\ell, 3} \le 3$. Also $\Disc_k = 1$ and $\zeta(2) L_{\ell / \Q}(3) \ge \zeta(2)^{1/2}$, so \eqref{eq:FieldBound} becomes
\begin{equation}\label{eq:FieldBoundQ}
\Disc_\ell \le \left( \frac{16 \pi^5}{\zeta(2)^{1/2} 288} \right)^{2/5} \approx 2.8116.
\end{equation}
The smallest discriminant of an imaginary quadratic field is $3$, so this rules out the case $k = \Q$.
\end{pf}

\begin{clm}
The totally real field $k$ must be $\Q(\al)$ and $\ell$ must be $k(\beta)$, where $\al^2 = 3$ and $\beta^2 = -1$.
\end{clm}

\begin{pf}
We know now that $n = 2$. As in previous cases, $h_{\ell, 3} \le 3$ and so the ideas used above give
\begin{equation}\label{eq:FieldBoundQuad}
\Disc_k^{1/2} \le \Disc_\ell^{1/4} \le \frac{(16 \pi^5)^{1/4}}{\zeta(4)^{1/16} 288^{1/8}} \approx 4.1011.
\end{equation}
Thus $\Disc_k \le 16$ and $\Disc_\ell \le 282$. There are $4$ such real quadratic fields, with discriminants $5$, $8$, $12$, and $13$, and nine such totally imaginary quartic fields \cite{Tables}. The cases where $\ell$ is a quadratic extension of $k$ are when $(\Disc_k, \Disc_\ell)$ is either $(5, 125)$, $(8, 256)$, or $(12, 144)$. We must rule out the first two cases.

For $(\Disc_k, \Disc_\ell) = (5, 125)$, we have that $k$ and $\ell$ both have class number one. We then use \eqref{eq:FieldBound} to see that
\begin{equation}\label{eq:FieldBoundn2d5}
\frac{\Disc_\ell^{5/2} \zeta_k(2) L_{\ell / k}(3)}{(16 \pi^5)^n \Disc_k h_{\ell, 3}} \ge \frac{125^{5/2} \zeta(4)^{1/2})}{5 (16 \pi^5)^2} \approx 0.00152 > \frac{1}{864},
\end{equation}
which rules out this case. The same calculation for $(8, 256)$ bounds the Euler characteristic from below by $0.00569$, which also rules out that case. This proves the claim.
\end{pf}

\begin{clm}
The algebra $D$ has degree $1$, i.e., $D = \ell$.
\end{clm}

\begin{pf}
Note that this is equivalent to showing that $\mc{T} = \mc{T}_0$ is empty. Let $\Gam < \PU(2, 1)$ be the maximal arithmetic lattice arising from the central simple $\ell$-algebra $D$ (with $k$ and $\ell$ as in the statement of the theorem). If $D \neq \ell$ (that is, if $D$ has degree $3$), then $\mc{T}_0$ is nonempty, i.e., $D$ must ramify above some finite nonempty set of places of $k$.

We then calculate
\[
e(\B^2 / \Gam) \ge \frac{1}{288} \prod_{v \in \mc{T}_0} e^{\prime \prime}(\para_v),
\]
with notation as above (see \cite[\S 8.2]{Prasad--Yeung} for the relevant special value of our $L$-function). To prove the claim, it then suffices to show that $e^{\prime \prime}(\para_v) > 1$ for all $v \in \mc{T}_0$. From the calculation of $e^{\prime \prime}(\para_v)$ in \cite[\S 2.5]{Prasad--Yeung}, for every $v \in \mc{T}_0$ we have $e^{\prime \prime}(\para_v) = 1$ if and only if the residue field of $k_v$ is $\F_2$, the field with $2$ elements.

Since $2$ ramifies in $k$, there is a unique place $v_2$ of $k$ over $2$ with residue field $\F_2$. This place is inert in $\ell$, but every place of $\mc{T}_0$ is necessarily split in $\ell$. (This means exactly that no central simple $\ell$ algebra of degree $3$ that ramifies above $2$ admits an involution of second kind.) Thus $e^{\prime \prime}(\para_v) > 1$ for all $v \in \mc{T}_0$ when it is nonempty, so $e(\B^2 / \Gam) > 1 / 288$. This proves the claim.
\end{pf}

Now we are only left with the commensurability class given in the statement of the theorem. We claim that there is a unique lattice of minimal covolume in this commensurability class. We first note that $\para$ must satisfy conditions (1)-(3) above, since $e^{\prime \prime}(\para_v) > 1$ whenever $\para_v$ does not satisfy all these assumptions. The claim now follows from the fact that no finite place in $k$ ramifies in $\ell$ as in \cite[$\S$9.6]{Prasad--Yeung}, which implies that there is a unique coherent collection $\para$ of parahoric subgroups of $G(\A_{k, f})$ that satisfy these conditions (i.e., nonuniqueness of $\para$ satisfying (1)-(3) can only arise from finite places of $k$ that ramify in $\ell$).

From \eqref{eq:IndexBound}, the fact that $\ell$ has class number one, and because we must have $\mc{T} = \emptyset$, we see that $[\wt{\Gam} : \Gam_{\para}]$ is either $1$ or $3$, hence the minimal possible orbifold Euler characteristic is either $1/288$ or $1/108$. One can prove that $[\wt{\Gam} : \Gam_{\para}] = 3$ algebraically, but it also follows immediately from the fact that there is a known lattice $\Gam$ in the commensurability class with $e(\B^2 / \Gam) = 1/288$, namely the Deligne--Mostow orbifold considered in this paper (e.g., see \cite[\S 7]{Sauter} or \S \ref{sec:DM} below). This completes the proof.
\end{pf}


\section{Deligne--Mostow orbifolds}\label{sec:DM}

\begin{figure}
\begin{center}
\includegraphics[width=30mm]{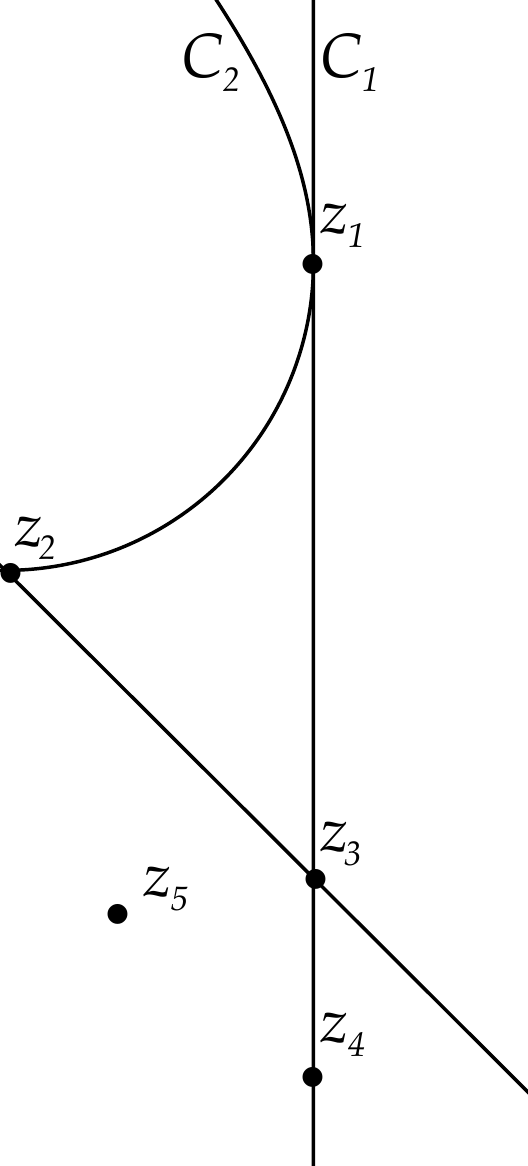}
\end{center}
\caption{The underlying space for $\mc{O}$}\label{fig:Orbifold}
\end{figure}

In this section, we give a geometric description of the Deligne--Mostow orbifold central to this paper. We describe it both as an analytic space and as a complex hyperbolic orbifold. Its existence (and known arithmeticity) suffice to finish the proofs of Theorems \ref{thm:DMVolume} and \ref{thm:ArithmeticAut}, but we give many more details for completeness. The treatment below (especially Figure \ref{fig:Orbifold}) is heavily influenced by conversations with Domingo Toledo, and any novel geometric insights should be considered his.

Recall that Hurwitz curves are exactly the manifold regular covers of the sphere with cone points of order $2$, $3$, and $7$, i.e., the $(2,3,7)$-triangle orbifold. Our replacement for the $(2, 3, 7)$-triangle orbifold is one of Deligne and Mostow's ball quotients. The Deligne--Mostow orbifolds are in many ways the natural generalization to $\mathbb B^2$ of the $(p, q, r)$-triangle orbifolds. In fact, in one complex dimension this construction exactly reproduces the triangle orbifolds \cite[\S 14.3]{Deligne--Mostow}.

We begin with an $(n+3)$-tuple of positive integers $(a_1, \dots, a_{n + 3})$ such that $\sum a_j = 2 t$ for some integer $t$. Set $\mu_j = a_j / t$. We then define the following condition on the $(n + 3)$-tuple $(\mu_1, \dots, \mu_{n + 3})$:
\begin{equation}
\mu_j + \mu_k < 1 \Rightarrow \left(1 - \mu_j - \mu_k \right)^{-1} \in \Z \tag{INT}
\end{equation}
Furthermore, we say that $(\mu_1, \dots, \mu_{n+3})$ satisfies ($\Sig$INT) if we weaken (INT) to the following:
\begin{eqnarray}
\mu_j + \mu_k < 1\ \textrm{and}\ \mu_j \ne \mu_k &\Rightarrow& \left(1 - \mu_j - \mu_k \right)^{-1} \in \Z \nonumber \\
\mu_j + \mu_k < 1\ \textrm{and}\ \mu_j = \mu_k &\Rightarrow& \left(1 - \mu_j - \mu_k \right)^{-1} \in \frac{1}{2} \Z \nonumber
\end{eqnarray}
If $(\mu_1, \dots, \mu_{n + 3})$ satisfies (INT) or ($\Sig$INT), we call it a \emph{ball tuple}.

For any ball tuple $(\mu_1, \dots, \mu_{n + 3})$, Deligne and Mostow \cite{Deligne--Mostow} and Mostow \cite{MostowSigma} constructed a lattice in $\PU(n, 1)$ with an explicit quotient orbifold described in terms of moduli of points on $\mathbb{P}^1$. For example, when $n = 4$ one can solve the equations
\[
r_j = \left(1 - \mu_j - \mu_{j + 1} \right)^{-1}
\]
for $j = 1, 2, 3$ to realize the $(r_1, r_2, r_3)$-triangle orbifold as the orbifold associated with the ball tuple $(\mu_1, \dots, \mu_4)$. See \cite{Kirwan--Lee--Weintraub} for an explicit description (in any dimension) of this orbifold as an algebraic space.

Our interest is in the orbifold $\mc{O}$ built from the ball tuple $(\frac{2}{12}, \frac{2}{12}, \frac{2}{12}, \frac{7}{12}, \frac{11}{12})$, which satisfies ($\Sig$INT) but not (INT). See \cite{Kirwan--Lee--Weintraub} for further details. Let $M$ be the space of distinct $5$-tuples on $\mathbb{P}^1$ and $Q$ the quotient of $M$ by the action of $\PGL_2$. An arbitrary $5$-tuple $(y_1, \dots, y_5)$ is called \emph{stable} (resp.~\emph{semistable}) if, for every $z \in \mathbb{P}^1$,
\[
\sum_{y_j = z} \mu_j < 1 \quad \textrm{(resp.~}\le 1\textrm{)}.
\]
If $M^{st}$ and $M^{sst}$ are the sets of stable and semistable $5$-tuples, respectively, we see (for our ball tuple) that $M \subset M^{st} = M^{sst}$. Since $\PGL_2$ also fixes the semistable points, we also obtain a quotient $Q^{sst}$ containing $Q$. It is a projective surface and one can show explicitly in coordinates that $Q^{sst}$ equals $\mathbb{P}^2$ (cf.~\cite[Thm.~4.1]{Kirwan--Lee--Weintraub}).

\begin{figure}
\begin{center}
\includegraphics[width=60mm]{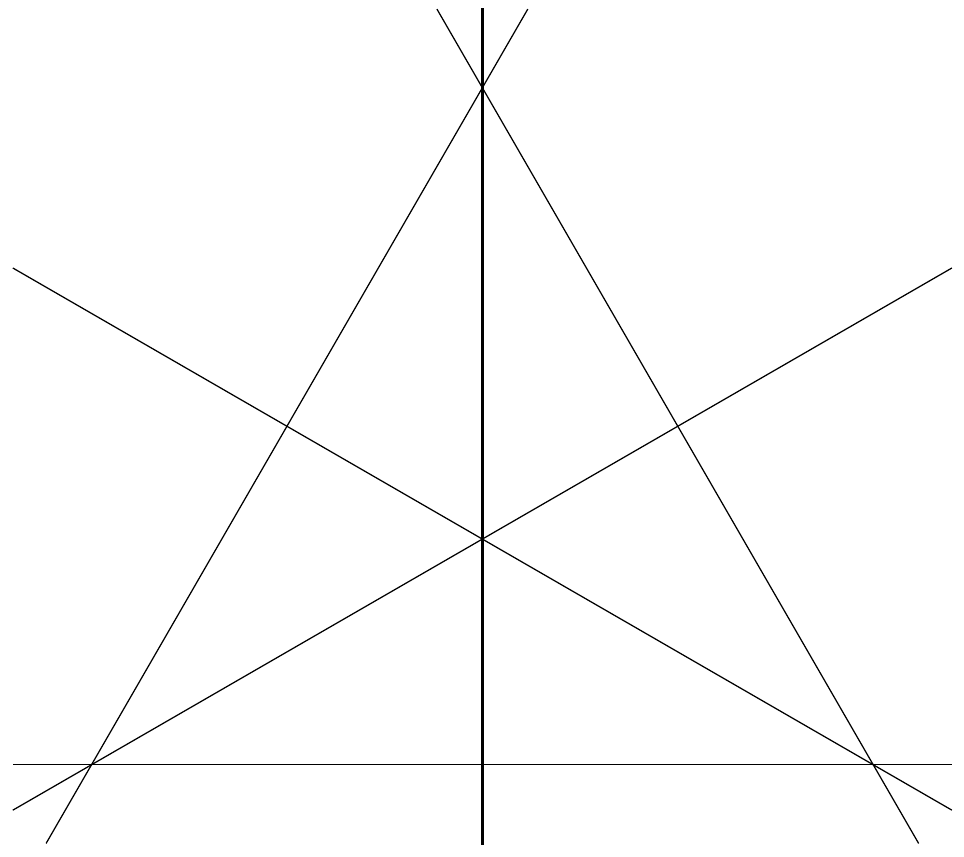}
\end{center}
\caption{Six lines on $\mathbb{P}^2$}\label{fig:SixLines}
\end{figure}

First, notice that identification of $Q^{sst}$ with $\mathbb{P}^2$ maps $Q^{sst} \ssm Q$ onto the classical arrangement of $6$ lines on $\mathbb{P}^2$ given in Figure \ref{fig:SixLines} (also see \cite{Hirzebruch}).
In particular, these lines are in one-to-one correspondence with pairs $(i, j)$ such that $\mu_i + \mu_j \le 1$, and the points of intersection between the lines correspond with triples $(i, j, k)$ such that $\mu_i + \mu_j + \mu_k \le 1$.

There is also an action of the symmetric group $\Sig_3$ on $Q^{sst}$, acting on the first $3$ coordinates of a $5$-tuple. This action clearly fixes $Q$, but also fixes $Q^{sst}$ since $\mu_1 = \mu_2 = \mu_3$. We can represent the quotient $Q^{sst}_\Sig$ of $Q^{sst}$ by this action as in Figure \ref{fig:Orbifold}. In coordinates, $C_1$ is the image of the three lines on $\mathbb{P}^2$ where $2$ coordinates are equal, and $C_2$ is the image of the three lines where one coordinate is zero.

Then $Q^{sst}_\Sig$ is the underlying space for a complex hyperbolic orbifold $\mc{O}$ with orbifold structure explicitly determined by the ball tuple. The orbifold locus is exactly the set of lines and points in Figure \ref{fig:Orbifold}, and we give the lifts to $\mathbb{P}^2$ of the points $z_j$ in Table \ref{tb:zjs}.
\begin{table}[h]
\begin{center}
\begin{tabular}{|c|c|c|}
\hline
$z_1$ & $[1 \,:\, 0 \,:\, 0],\, [0 \,:\, 1 \,:\, 0],\, [0 \,:\, 0 \,:\, 1]$ \\
\hline
$z_2$ & $[1 \,:\, 1 \,:\, 1]$ \\
\hline
$z_3$ & $[1 \,:\, 1 \,:\, 0],\, [1 \,:\, 0 \,:\, 1],\, [0 \,:\, 1 \,:\, 1]$ \\
\hline
$z_4$ & $[0 \,:\, 1 \,:\, -1],\,[1 \,:\, 0 \,:\, -1],\,[1 \,:\, -1 \,:\, 0]$ \\
\hline
$z_5$ & $[1\,:\,\om\,:\,\om^2],\, [1\,:\,\om^2\,:\,\om]$ where $\om^3 = 1$ \\
\hline
\end{tabular}
\caption{Lifts of the points $z_j$}\label{tb:zjs}
\end{center}
\end{table}
Orbifold weights of the curves $C_1, C_2$ can be directly computed via the above correspondence between the singular locus and ($\Sig_3$-equivalence classes of) pairs $(i, j)$ with $\mu_i + \mu_j < 1$. Since our ball tuple satisfies ($\Sig$INT) but not (INT), the curve $C_1$ has orbifold weight $4 = (1 - (2 + 7)/12)^{-1}$ and $C_2$ has weight $3 = 2(1 - (2 + 2)/12)^{-1}$. The local group for each of these curves is cyclic of the appropriate order, generated by a complex reflection. The marked points have local groups generated by reflections, and one can identify the local groups from calculating the associated Gram matrix in appropriate coordinates and consulting \cite[Ch.~6]{Lehrer--Taylor}. The orbifold weight and local group for each $z_i$ is given in Table \ref{tb:LocalGps}.
\begin{table}[h]
\begin{center}
\begin{tabular}{|c|c|c|}
\hline
Point & Orbifold weight & Local Group \\
\hline
$z_1$ & $288$ & $G_{10}$ \\
$z_2$ & $24$ & $G_{4}$ \\
$z_3$ & $12$ & $\Z / 3 \Z \times \Z / 4 \Z$ \\
$z_4$ & $8$ & $\Z / 4 \Z \times \Z / 2 \Z$ \\
$z_5$ & $3$ & $\Z / 3 \Z$ \\
\hline
\end{tabular}
\caption{Local groups at orbifold points}\label{tb:LocalGps}
\end{center}
\end{table}
Here $G_m$ denotes the $m^{th}$ exceptional group in the Shephard--Todd classification of finite complex reflection groups.

One can compute from this information that the orbifold Euler characteristic of $\mc{O}$ is $1 / 288$. Also, see the tables in \cite{Sauter}, where one sees that this ball tuple gives the Deligne--Mostow orbifold of smallest orbifold Euler characteristic. We also know from \cite{MostowSigma} that the orbifold fundamental group of $\mc{O}$ determines an arithmetic lattice in $\PU(2, 1)$ commensurable with the lattice $\Gam$ in Theorem \ref{thm:MinVol}. By uniqueness of the arithmetic lattice of minimal covolume, $\pi_1^{\textrm{orb}}(\mc{O})$ must equal $\Gam$. This completes the proof of Theorems \ref{thm:DMVolume} and \ref{thm:ArithmeticAut}. \qed


\begin{rems}\ 
\begin{enumerate}

\item The tables in \cite{Sauter} show that there is no known nonarithmetic counterexample to Conjecture \ref{MinVolConjecture} amongst the Deligne--Mostow lattices. Martin Deraux informed me that none of the new nonarithmetic lattices of Deraux--Parker--Paupert \cite{Deraux--Parker--Paupert} are counterexamples either.

\item We also see from the tables in \cite{Sauter} that this orbifold is one of the orbifolds constructed by Mostow using complex reflection groups \cite{MostowReflect}. In Mostow's notation, the group is $\Gam(3, \frac{1}{3})$. Since these groups are often called complex hyperbolic triangle groups, this also makes the orbifold in Theorems \ref{thm:ArithmeticAut} and \ref{thm:DMVolume} a natural generalization of the $(2, 3, 7)$-triangle orbifold.

\end{enumerate}
\end{rems}


\section{Extremal automorphism groups}\label{sec:Gps}


We recall again that a Hurwitz group is a finite group $F$ of order $84(g - 1)$ for some $g \ge 2$ such that $F = \Aut(C)$ for some smooth projective curve $C$ of genus $g$. Hurwitz groups are precisely the finite quotients of the $(2, 3, 7)$-triangle group with torsion-free kernel. Such groups are in one-to-one correspondence with finite groups generated by two elements $x, y$ such that $x$ has order $2$, $y$ has order $3$ and $xy$ has order $7$. See \cite{Conder} for more on these groups.

Using similar reasoning, we can classify the groups that appear as symmetry groups of arithmetic ball quotients $S$ with $|\Aut(S)| = 288 e(S)$. Let $\Gam < \PU(2, 1)$ be the associated maximal arithmetic lattice, i.e., the orbifold fundamental group of the Deligne--Mostow orbifold central to this paper. One can, for example, classify these automorphism groups via the presentation of Cartwright and Steger for $\Gam$ \cite{Cartwright--Steger} or Mostow's presentation of the lattice as a complex reflection group \cite{MostowReflect}.

We use Cartwright and Steger's presentation:
\begin{eqnarray}
\Gam = \big\langle b,\, j,\, u,\, v ~:~ &&b^3,\ u^4,\ v^8,\ (u, j),\ (v, j),\ j^{-3} v^2, \nonumber \\
 &&(b v u^3)^3,\ u v (u v^{-1})^2,\ (b j)^2 u^{-2} v^{-1},\ b^{-1} u^{-2} v^{-1} b v u^2 \big\rangle, \nonumber
\end{eqnarray}
where $(x,y)$ denotes the commutator of $x$ and $y$. MAGMA \cite{MAGMA} code (based on the code of Cartwright and Steger) implementing what follows is available on the author's webpage. What is nice about this presentation from our perspective is that one can show that the subgroup of $\Gam$ generated by $j$, $u$, and $v$ is isomorphic to the group $G_{10}$. Recall that $G_{10}$ is the exceptional complex reflection group of order $288$ in the Shephard--Todd classification, and is the local group of the point $z_1$ in Figure \ref{fig:Orbifold}. Thus $\Gam$ is generated by $G_{10}$ and an element of order $3$, which gives the following analogue to the group theoretic classification of Hurwitz groups.


\begin{prop}\label{prop:GpCharacterization}
Automorphism groups of smooth arithmetic ball quotients $S$ with $|\Aut(S)| = 288 e(S)$ are in one-to-one correspondence with finite groups $H$ such that:
\begin{enum}

\item $H$ is generated by four elements, $\conj{b}$, $\conj{j}$, $\conj{u}$, $\conj{v}$ such that $\conj{b}$ has order $3$ and the subgroup of $H$ generated by $\conj{j}$, $\conj{u}$, and $\conj{v}$ is isomorphic to the complex reflection group $G_{10}$;

\item $(\conj{b}\, \conj{v}\, \conj{u}^3)^3 = (\conj{b}\, \conj{j})^2 \conj{b}\, \conj{j}\, \conj{u}^{-2} \conj{v}^{-1} = \conj{b}^{-1} \conj{u}^{-2} \conj{v}^{-1} \conj{b}\, \conj{v}\, \conj{u}^2 = 1$;

\item if $w(b, j, u, v)$ is a word representing a conjugacy class of elements in $\Gam$ of finite order $n > 1$, then $w(\conj{b}, \conj{j}, \conj{u}, \conj{v})$ has order $n$ in $H$.

\end{enum}
\end{prop}


The list of conjugacy classes of finite order elements of $\Gam$ is a bit long, so we omit the list from the statement of the proposition. Equivalently, one must consider all the finite groups in Table \ref{tb:LocalGps}. For the interested reader, the author's MAGMA code enumerates these elements, as does Cartwright and Steger's.


\begin{pf}[Proof of Proposition \ref{prop:GpCharacterization}]
If $S$ is a smooth arithmetic ball quotient for which $|\Aut(S)| = 288e(S)$, then
\[
e(S / \Aut(S)) = e(S) / |\Aut(S)| = 1 / 288.
\]
By Theorem \ref{thm:DMVolume}, the orbifold fundamental group of $S / \Aut(S)$ must be $\Gam$, hence there is a surjective homomorphism $\rho : \Gam \to H$ with kernel isomorphic to $\pi_1(S)$. Then $\conj{b} = \rho(b)$ etc.~is a generating set for $H$ satisfying (1) and (2) in the statement of the proposition. Since $\pi_1(S)$ is torsion-free, every element of finite order $n$ in $\Gam$ must map onto an element of order $n$ in $H$, which means that (3) holds.

Conversely, given a group $H$ satisfying (1)-(3), this exactly means that there is a surjective homomorphism $\Gam \to H$ with torsion-free kernel $\Lam$. Let $S$ be the associated smooth ball quotient $\B^2 / \Lam$. Then $S$ admits an action by $H$ with $S / H = \B^2 / \Gam$, so $H$ is isomorphic to a subgroup of $\Aut(S)$. It follows that
\[
\Aut(S) \ge |H| = e(S) / e(\B^2 / \Gam) = 288 e(S),
\]
and equality follows from Theorem \ref{thm:ArithmeticAut}.
\end{pf}


Though many groups are the automorphism groups of Hurwitz curves \cite{Conder}, the first several Hurwitz curves come from congruence subgroups of $\Delta_{2, 3, 7}$, including the Klein quartic, which is the smallest genus Hurwitz curve, and the Macbeath curve, which has automorphism group $\PSL_2(\F_8)$ and genus $7$. A similar result, which we now explain, holds for ball quotients.

Using MAGMA, we computed all the smooth arithmetic ball quotients $S$ with $|\Aut(S)| = 288 e(S)$ and $e(S) \leq 375$. This was done by enumerating the finite index normal subgroups of $\Gam$ with index at most $108000$, then checking torsion-freeness by checking that no conjugacy class of finite order elements of $\Gam$ is contained in a given subgroup.

There are exactly two such surfaces, which we record in Table \ref{tb:First3} along with $H_1(S, \Z)$, the irregularity $q$, geometric genus $p_g$, and $h^{1,1}$. We computed $H_1$ using MAGMA by finding the abelianization of the associated subgroup of $\Gam$, and the other invariants can be computed easily from there. Recall from above that $K_S^2 = 3 e(S)$ and $\chi(\mc{O}_S) = e(S) / 3$ as in Corollary \ref{cor:BallInvariantAut}.
\begin{table}[h]
\begin{center}
\begin{tabular}{|c|c|c|c|c|c|}
\hline
 & $e$ & $H_1(S, \Z)$ & $q$ & $p_g$ & $h^{1,1}$ \\
\hline
$S_1$ & $63$ & $\Z^{14}$ & $7$ & $27$ & $35$ \\
\hline
$S_2$ & $252$ & $(\Z / 2 \Z)^{15} \oplus \Z^{14}$ & $7$ & $90$ & $98$ \\
\hline
\end{tabular}
\caption{The two smallest arithmetic $S$ with $|\Aut(S)| = 288 e(S)$}\label{tb:First3}
\end{center}
\end{table}
They are the smallest arithmetic ball quotients satisfying $|\Aut(S)| = 288 e(S)$. If Conjecture \ref{AutConjecture} holds, these are the two Hurwitz ball quotients of smallest Euler characteristic.

We now prove Theorem \ref{thm:Smallest}, which identifies $S_1$ as a quotient of the ball by a principal congruence lattice (hence as a PEL Shimura variety) with automorphism group $\PSU(3, \F_3) \times (\Z / 3 \Z)$.


\begin{pf}[Proof of Theorem \ref{thm:Smallest}]
In order to describe $S_1$ more precisely, we need notation from \S \ref{sec:Arithmetic}. Consider:
\begin{enum}

\item $k = \Q(\al)$, where $\al^2 = 3$

\item $\ell = k(\beta)$, where $\beta^2 = -1$

\item $\tau$ the nontrivial Galois involution of the quadratic extension $\ell / k$

\item $h$ the $\tau$-hermitian form on $\ell^3$ with matrix
\[
\begin{pmatrix}
1 & 0 & 0 \\ 0 & 1 & 0 \\ 0 & 0 & 1 - \al
\end{pmatrix}
\]

\end{enum}
Recall that this data uniquely determines the absolutely almost simple and simply connected algebraic group $G$ associated with $\Gam$, where $\Gam$ is the orbifold fundamental group of the Deligne--Mostow orbifold with ball tuple $(\frac{2}{12}, \frac{2}{12}, \frac{2}{12}, \frac{7}{12}, \frac{11}{12})$.

Furthermore, let $\para$ be the unique coherent collection of parahoric subgroups of $G(\mathbb{A}_{k, f})$ satisfying (1)-(3) in the proof of Theorem \ref{thm:MinVol}. Recall that uniqueness follows from the fact that $\ell / k$ is unramified at all finite places. Let $\Gam_{\para}$ be the associated arithmetic subgroup of $G(k)$ and $\conj{\Gam}_{\para}$ its image in $\PU(2, 1)$.

The proof of Theorem \ref{thm:MinVol} implies that $\Gam$ is the normalizer of $\conj{\Gam}_{\para}$ in $\PU(2, 1)$, and
\[
[\Gam : \conj{\Gam}_{\para}] = 3.
\]
In more classical terms, one can realize $\conj{\Gam}_{\para}$ as
\[
\PSU(h, \mc{O}_k) = \left\{ x \in \PSL_3(\mc{O}_\ell)~:~{}^t\tau(x) h x = h \right\},
\]
where $\mc{O}_F$ is the ring of integers in the field $F$ and ${}^t$ denotes matrix transpose. Briefly, for each each finite place $v$ of $k$, $\para_v$ is the stabilizer in the associated Bruhat--Tits building of the lattice $\mc{L}_v$ generated by the lattice $\mc{O}_\ell^3 \subset \ell^3$. Thus $\Gam_{\para}$ is the special unitary group of the hermitian lattice $(\mc{O}_\ell^3, h)$. Using MAGMA, one sees that $\Gam$ has a unique normal subgroup of index $3$, which then must be $\conj{\Gam}_{\para}$.

There is a unique prime ideal $\mf{p}_3$ of $\mc{O}_k$ over $3$, and this prime is inert in $\ell$. Therefore, we have a reduction homomorphism
\[
r_{\mf{p}_3} : \conj{\Gam}_{\para} \to \PSU(3, \F_3),
\]
where $\PSU(3, \F_3)$ is the unique projective special unitary group in $3$ variables for the extension of finite fields $\F_9 / \F_3$. Let $\conj{\Gam}(\mf{p}_3)$ be the kernel of $r_{\mf{p}_3}$. Then $\conj{\Gam}(\mf{p}_3)$ is a principal congruence subgroup of $\conj{\Gam}_{\para}$ in the usual sense. We can also realize $\conj{\Gam}(\mf{p}_3)$ as $\conj{\Gam}_{\para^{(3)}}$ for $\para^{(3)} \subset \para$ a proper open compact subgroup. It is well-known that $\conj{\Gam}(\mf{p}_3)$ is torsion-free, and a simple calculation shows that
\[
e(\B^2 / \conj{\Gam}(\mf{p}_3)) = 63.
\]
To prove that $\B^2 / \conj{\Gam}(\mf{p}_3) = S_1$, we must show that $\conj{\Gam}(\mf{p}_3)$ is a normal subgroup of $\Gam$.

Using MAGMA, one sees that there is (up to automorphisms) a unique homomorphism from $\conj{\Gam}$ to $\PSU(3, \F_3)$, which allows us to precisely identify $\conj{\Gam}(\mf{p}_3)$. From this, one can check that $\conj{\Gam}(\mf{p}_3)$ is a normal subgroup of $\Gam$. This identifies $S_1$ as the quotient of the ball by a principal congruence subgroup. To give the precise description of $S_1$ as a moduli space of principally polarized abelian varieties with $\ell$-endomorphisms and $\para^{(3)}$ level structure, we refer the reader to \cite[Thm.~8.17]{Milne}. It remains to identify $\Gam / \conj{\Gam}(\mf{p}_3)$ with $\PSU(3, \F_3) \times (\Z / 3 \Z)$, which one can do directly by having MAGMA construct an isomorphism between them.
\end{pf}


Using further MAGMA calculations like those in the proof of Theorem \ref{thm:Smallest}, one can show that the surface $S_2$ is a $4$-fold regular cover of $S_1$. Moreover, one can check that its automorphism group is isomorphic to $\PSU(3, \F_3) \times \mathrm{A}_4$, where $\mathrm{A}_4$ is the alternating group on four letters. We were not able to show that $\pi_1(S_2)$ contains a principal congruence lattice, so it may be a noncongruence subgroup of $\conj{\Gam}_{\para}$. We close with two final remarks.


\begin{rems}\ 
\begin{enum}

\item These surfaces are also commensurable with the smooth compact arithmetic ball quotient $S$ with $e(S) = 3$ and first Betti number $2$ discovered by Cartwright and Steger in their classification of fake projective planes \cite{Cartwright--Steger}. Their surface does not have the extremal symmetry properties of the above surfaces ($\Aut(S) \cong \Z / 3 \Z$), but it is a quotient of $S_1$ by an action of the Frobenius group of order $21$. In fact, one can verify with MAGMA that the covering group is a subgroup of $\PU(3, \F_3)$, that is, $S$ is an intermediate covering between $S_1$ and $\B^2 / \conj{\Gam}_{\para}$ (with notation as in the proof of Theorem \ref{thm:Smallest}). Consequently, there should be a modular interpretation of the Cartwright--Steger surface, which we plan to explore in joint work with Domingo Toledo.

\item We know very little about the problems considered in this paper for higher-dimensional ball quotients. We do make one remark, which follows from a combination of work of the author with Vincent Emery \cite{Emery--Stover} along with recent work of McMullen \cite{McMullen}. Let $\mc{O} = \B^n / \Lam$ for some nonuniform arithmetic lattice $\Lam < \PU(n, 1)$ (for any $n$). Then
\[
|e(\mc{O})| \ge \frac{809}{5746705367040}
\]
with equality if and only $\mc{O}$ is the unique $9$-dimensional Deligne--Mostow orbifold (see also \cite{Allcock}).
\end{enum}
\end{rems}


\bibliography{hurwitz}


\end{document}